\documentclass[a4paper]{amsproc}

\usepackage{amsmath}
\usepackage{amsfonts}
\usepackage{amssymb}
\usepackage{pstricks}
\usepackage{graphicx}
\usepackage{latexsym}
\usepackage{geometry}
\usepackage{url}

\theoremstyle{plain}

\theoremstyle{definition}

\numberwithin{equation}{section}
\renewcommand{\leq}{\leqslant}
\renewcommand{\geq}{\geqslant}

\newcommand{\atan}{\operatorname{atan}}

\title[\hspace{-10pt}Corrections, Improvements, Comments on \textit{Gradshteyn \& Ryzhik} Integrals]{Corrections, Improvements, and Comments on Some \textit{Gradshteyn and Ryzhik} Integrals}



\author[Robert C. Elliott]{Robert C. Elliott${}^{*}$}

\author{Witold A. Krzymie\'{n}}

\address{\hfill{\it Received 01 09
2023, revised 30 11 2023}\newline Department of Electrical and Computer Engineering, University of Alberta, \newline 11th Floor, D-ICE Bldg., 9211–116th St. NW, Edmonton, AB  T6G~1H9, Canada}
\email{rce@ualberta.ca (${}^{*}$Corresponding author), krzymien@ualberta.ca}

\begin{document}

{\begin{flushleft}\baselineskip9pt\scriptsize {\bf SCIENTIA}\newline
Series A: {\it Mathematical Sciences}, Vol. 34 (2024), xx--xx
\newline Universidad T\'ecnica Federico Santa Mar{\'\i}a
\newline Valpara{\'\i}so, Chile
\newline ISSN 0716-8446
\newline {\copyright\space Universidad T\'ecnica Federico Santa
Mar{\'\i}a\space 2024}
\end{flushleft}}

\vspace{10mm} \setcounter{page}{1} \thispagestyle{empty}

\begin{abstract}
In this paper, we prove that two integrals from Gradshteyn and Ryzhik (2014) \cite{gr14} (namely, Eqs. 3.937 1 and 3.937 2) provide incorrect results in certain conditions. We derive those conditions herein and provide the corrections required for those two formulas. We furthermore derive improved formulas for the solutions to those integrals that are less complicated, avoid the errors of the original formulas, and work under a larger range of parameter values. The improved formulas are used to verify the results of a few other related integrals from \cite{gr14}; the previous results need correction in some instances, or are correct but can be extended in other instances. Lastly, we also consider the extended case of complex-valued parameters and derive the resulting formulas.
\end{abstract}

\maketitle

\section{Introduction}\label{sec:intro}
In Gradshteyn and Ryzhik's \textit{Table of Integrals, Series, and Products} \cite{gr14}, the following results for two definite integrals (\cite[Eq.~3.937~1]{gr14} and \cite[Eq.~3.937~2]{gr14}, respectively) are listed:
\begin{equation}\label{eq:sin_int_old}
\begin{aligned}
&\int_{0}^{2\pi} \exp(p \cos x + q \sin x) \sin(a \cos x + b \sin x - mx) \, dx \\
&\quad= i\pi\big[(b-p)^2 + (a+q)^2\big]^{-m/2}  \\
&\qquad \times \left[(A+iB)^{m/2\,} I_m \!\left(\!\sqrt{C-iD}\right) - (A-iB)^{m/2\,} I_m \!\left(\!\sqrt{C+iD}\right)\right].
\end{aligned}
\end{equation}
\begin{equation}\label{eq:cos_int_old}
\begin{aligned}
&\int_{0}^{2\pi} \exp(p \cos x + q \sin x) \cos(a \cos x + b \sin x - mx) \, dx  \\
&\quad= \pi\big[(b-p)^2 + (a+q)^2\big]^{-m/2}  \\
&\qquad \times \left[(A+iB)^{m/2\,} I_m \!\left(\!\sqrt{C-iD}\right) + (A-iB)^{m/2\,} I_m \!\left(\!\sqrt{C+iD}\right)\right].
\end{aligned}
\end{equation}
$I_m(z)$ is the modified Bessel function of the first kind. For these two equations, it is specified that $m \in \mathbb{N}$ (i.e., $m = 0, 1, 2, \dots$), $A = p^2 - q^2 + a^2 - b^2$, $B = 2(pq + ab)$, $C = p^2 + q^2 - a ^2 - b^2$, $D = 2(ap + bq)$, and it is required that $(b-p)^2 + (a+q)^2 > 0$. The original source for these two formulas is given as Gr\"{o}bner and Hofreiter (\cite[Sec.~339, Eq.~9b]{gh73} and \cite[Sec.~339, Eq.~9a]{gh73}, respectively). Note that $D$ as given here corrects a typo in the 8th edition of Gradshteyn and Ryzhik \cite{gr14}, which accidentally inserted a minus sign into the equation for $D$. The correct version without the minus sign is in Gr\"{o}bner and Hofreiter \cite{gh73} and earlier editions of \cite{gr14}.

In working with these two equations, we have determined that both of them can yield an incorrect result. Specifically, when $m$ has an odd value, the sign of the result will be the opposite of what it should be. The conditions for this occurring depend on what the values of $a$, $b$, $p$, and $q$ are.
In what follows, we show the source of the error and derive the conditions in which the sign error occurs. We shall also derive improved alternative formulas that are shorter, avoid this error, and do not require $(b-p)^2 + (a+q)^2 > 0$.

\section{The Original Derivation}\label{sec:orig_dev}
Gr\"{o}bner and Hofreiter \cite{gh73} give a hint about how to obtain the original formulas, namely, by contour integration. Let us combine (\ref{eq:cos_int_old}) and (\ref{eq:sin_int_old}) together into one equation by using $e^{i\theta} = \cos \theta + i \sin \theta$. This gives:
\begin{equation}\label{eq:orig_int}
\begin{aligned}
    f &= \int_{0}^{2\pi} \exp(p \cos x + q \sin x) \exp[\,i(a \cos x + b \sin x - mx)]\, dx \\
      &= \int_{0}^{2\pi} \exp\!\big[ (p+ia)\cos x +(q+ib)\sin x\big] e^{-imx} \, dx.
\end{aligned}
\end{equation}
The results of the integral in (\ref{eq:cos_int_old}) will then ultimately be given by $\mathfrak{Re}(f) = (f+\overline{f})/2$, while the results of the integral in (\ref{eq:sin_int_old}) will be given by $\mathfrak{Im}(f) = (f-\overline{f})/(2i)$, where $\overline{f}$ denotes the complex conjugate of $f$.

We next make the substitution $z = e^{ix}$, which transforms the integral into a contour integral with the contour $\mathcal{C}$ being the complex unit circle.
\begin{equation}\label{eq:cont_int_pt1}
\begin{aligned}
    f &= \oint_{\mathcal{C}} \exp\!\left[ (p+ia)\left(\frac{z+z^{-1}}{2}\right) +(q+ib)\left(\frac{z-z^{-1}}{2i}\right)\right] z^{-m}\, \frac{dz}{iz} \\
    &= \oint_{\mathcal{C}} \exp\!\left[ \left(\frac{p+ia}{2}\right)(z+z^{-1}) + \left(\frac{b-iq}{2}\right)(z-z^{-1})\right] z^{-m-1}\, \frac{dz}{i}\\
    &= \frac{1}{i}\oint_{\mathcal{C}} \exp\!\left[ \frac{(p+b)+i(a-q)}{2}\, z + \frac{(p-b)+i(a+q)}{2}\, z^{-1}\right] z^{-m-1}\, dz
\end{aligned}
\end{equation}
For brevity, we shall denote $X = \big[(p{+}b)+i(a{-}q)\big]/2$ and $Y = \big[(p{-}b)+i(a{+}q)\big]/2$.

Next, we substitute in the power series representation $e^x = \sum_{n=0}^{\infty} (x^n/n!)$:
\begin{equation}\label{eq:cont_int_pt2}
\begin{aligned}
    f &= \frac{1}{i}\oint_{\mathcal{C}}\, \sum_{n=0}^{\infty} \frac{z^{-m-1} (Xz+Yz^{-1})^n}{n!} \, dz \\
    &= \frac{1}{i}\oint_{\mathcal{C}}\, \sum_{n=0}^{\infty} \frac{z^{-m-1}}{n!} \left[ \sum_{\ell = 0}^{n} \binom{n}{\ell}(Xz)^{\ell} (Yz^{-1})^{n-\ell}\right] dz \\
    &= \frac{1}{i}\oint_{\mathcal{C}}\, \sum_{n=0}^{\infty} \sum_{\ell = 0}^{n} \frac{z^{2\ell -m-n-1}}{\ell!\, (n-\ell)!} X^{\ell}\, Y^{n-\ell}\, dz \\
    &\overset{\textrm{(a)}}{=} \frac{1}{i}\oint_{\mathcal{C}}\, \sum_{\ell=0}^{\infty} \sum_{n = \ell}^{\infty} \frac{z^{2\ell -m-n-1}}{\ell!\, (n-\ell)!} X^{\ell}\, Y^{n-\ell}\, dz \\
    &\overset{\textrm{(b)}}{=} \frac{1}{i}\oint_{\mathcal{C}}\, \sum_{\ell=0}^{\infty} \sum_{k = -\infty}^{\ell} \frac{z^{k-m-1}}{\ell!\, (\ell - k)!} X^{\ell}\, Y^{\ell - k}\, dz  \\
    &\overset{\textrm{(c)}}{=} \frac{1}{i}\oint_{\mathcal{C}}\, \sum_{k = -\infty}^{\infty} \left[\sum_{\ell = \max(0,k)}^{\infty} \frac{X^{\ell}\, Y^{\ell - k}}{\ell!\, (\ell - k)!} \right]z^{k-m-1}\, dz.
    \end{aligned}
\end{equation}
Equalities (a) and (c) come from reversing the order of the summations, while equality (b) substitutes $k=2\ell - n$. We are left with a contour integral over a Laurent series of $z$. From the residue theorem for contour integration, the value of the integral will be $2\pi i$ times the coefficient of the $z^{-1}$ term of the series; see e.g. Krantz \cite{krantz99}. This coefficient occurs when $k=m$. We therefore obtain:
\begin{equation}\label{eq:cont_int_result}
    f = \frac{2\pi i}{i} \sum_{\ell = m}^{\infty} \frac{X^{\ell}\, Y^{\ell - m}}{\ell!\, (\ell - m)!} = 2\pi \sum_{j = 0}^{\infty} \frac{X^{j+m}\, Y^{j}}{(j+m)!\, j!} = 2\pi X^m \sum_{j = 0}^{\infty} \frac{X^{j}\, Y^{j}}{j!\,(j+m)!},
\end{equation}
with the substitution $j = \ell - m$.

We now note the following power series representation of $I_{\nu}(z)$ (Gradshteyn and Ryzhik \cite[Eq.~8.445]{gr14}):
\begin{equation}\label{eq:Iv_def}
    I_{\nu}(z) = \sum_{k=0}^{\infty} \frac{(z/2)^{\nu + 2k}}{k!\, \Gamma(\nu + k + 1)} = \left(\frac{z}{2}\right)^{\nu} \sum_{k=0}^{\infty} \frac{(z/2)^{2k}}{k!\, \Gamma(\nu + k + 1)}.
\end{equation}
Note also that for $\nu \in \mathbb{N}$, $\Gamma(\nu + k + 1) = (\nu + k)!$. Then, a comparison of (\ref{eq:cont_int_result}) and (\ref{eq:Iv_def}) shows they are very similar in form. We may manipulate (\ref{eq:cont_int_result}) to obtain
\begin{equation}\label{eq:f_manip}
\begin{aligned}
    f &= 2\pi (X^{1/2})^{m\,} (\overline{Y}^{1/2})^{m\,} (\overline{Y}^{1/2})^{-m\,} (Y^{1/2})^{-m\,}\\[-5pt]  &\hspace{100pt} \times (Y^{1/2})^{m\,} (X^{1/2})^{m\,} \sum_{j = 0}^{\infty} \frac{(X^{1/2})^{2j}\, (Y^{1/2})^{2j}}{j!\,(j+m)!}.
\end{aligned}
\end{equation}
Then, combining alike exponents:
\begin{equation}\label{eq:f_error_combine}
    f = 2\pi (X\overline{Y})^{m/2} (Y\overline{Y})^{-m/2\,}  \big[(XY)^{1/2}\big]^{m} \sum_{j = 0}^{\infty} \frac{\big[(XY)^{1/2}\big]^{2j}}{j!\,(j+m)!},
\end{equation}
and by assigning $z/2 = (XY)^{1/2}$, we get
\begin{equation}\label{eq:f_Im_sub}
    f = 2\pi (X\overline{Y})^{m/2} (Y\overline{Y})^{-m/2\,}  I_m\big(2\sqrt{XY}\big).
\end{equation}

From basic algebra, one can find that $Y\overline{Y} = |Y|^2 = [(b{-}p)^2 + (a{+}q)^2]/4$, $X\overline{Y} = (A-iB)/4$, and $2\sqrt{XY} = \sqrt{C+iD}$. The factor of 1/4 cancels in the first two terms of (\ref{eq:f_Im_sub}). We then finally obtain
\begin{equation}\label{eq:f_Im_final}
    f = 2\pi \big[(b-p)^2 + (a+q)^2\big]^{-m/2}\, (A-iB)^{m/2} I_m\big(\sqrt{C+iD}\big).
\end{equation}
It then just remains to use the properties $\overline{(zw)} = \overline{z}\,\overline{w}$, $\overline{z^{\alpha}} = \overline{z}^{\hspace{1pt}\alpha}$ for $\alpha \in \mathbb{R}$, and $I_{\nu}(\overline{z}) = \overline{I_{\nu}(z)}$ for $\nu \in \mathbb{R}$ (NIST \cite[Eq.~10.34.7]{NIST_DLMF}). The result in (\ref{eq:cos_int_old}) is obtained by taking $(f+\overline{f})/2$, and the result in (\ref{eq:sin_int_old}) is obtained by taking $(f-\overline{f})/(2i)$. Due to the factor of $Y^{-m/2}$ in (\ref{eq:f_manip}), which propagates to later equations, the formulas cannot be used if $Y=0$ (or equivalently $|Y|^2=0$) to avoid dividing by zero. The solutions therefore require $(b-p)^2 + (a+q)^2 > 0$.

\section{The Error}\label{sec:error}
The reader may have intuited the source of the error in the derivation by this point. Perhaps surprisingly, the error is not itself in breaking apart terms into products of square roots in (\ref{eq:f_manip}), though the act of doing so is eventually what leads to the error. Rather, the error is in the recombination of terms in (\ref{eq:f_error_combine}). In general, raising a complex number to a power is not well-defined, as the result may be multi-valued. (This is experienced even with purely real values. For instance, it is well known that $x^{1/2}$ in fact has two valid solutions: the positive and negative square roots.) A unique and/or principal value for $z^{\alpha}$ can be defined in two cases: 1) if $z$ is complex, $\alpha$ is an integer, and $z \neq 0$ if $\alpha \leq 0$, then $z^{\alpha}$ is single-valued; 2) if $z$ is positive real and $\alpha$ is complex, then $z^{\alpha} = e^{\alpha \log z}$. Unfortunately, neither of these cases may apply in general for this problem.

The result of this is that $z^{\alpha} w^{\alpha}$ does not always equal $(zw)^{\alpha}$. A simple counterexample is $z{\,=\,}w{\,=\,}{-1}$ and $\alpha = 1/2$. In this case, $z^{\alpha} = w^{\alpha} = i$, so $z^{\alpha}w^{\alpha} = i{\,\times\,}i = -1$. However, $(zw)^{\alpha} = ({-1}{\,\times\,}{-1})^{1/2} = 1^{1/2} = 1$. One may denote $z$ as $|z|e^{i \theta_z}$, where $\theta_z = \arg z$. One could therefore define $z^{\alpha}$ as $|z|^{\alpha} e^{i \alpha \theta_z}$, and thus $z^{\alpha}w^{\alpha} = (|z|{\cdot}|w|)^{\alpha} e^{i \alpha (\theta_z + \theta_w)}$. However, $z$ may be equally be denoted as $|z|e^{i (\theta_z + 2n\pi)}$ for any integer $n$. $\theta_z$, being the \textit{principal} argument of $z$, corresponds to the branch where $n = 0$, such that $\theta_z \in (-\pi,\pi]$. One may also denote $y = (zw)^{\alpha} = |y| e^{i \theta_y} = (|z|{\cdot}|w|)^{\alpha} e^{i \alpha \theta_y}$. It then happens that $z^{\alpha} w^{\alpha}$ may not equal $(zw)^{\alpha}$ if $\alpha (\theta_z + \theta_w)$ crosses over the branch cut at $\pi$ radians and ends up in a different branch than $\alpha \theta_y$. Specifically, if $\alpha$ is an odd integer multiple of $1/2$, then $z^{\alpha} w^{\alpha} \neq (zw)^{\alpha}$ if $\theta_z + \theta_w \not\in (-\pi,\pi]$, or in other words, if either $\theta_z + \theta_w > \pi$ or $\theta_z + \theta_w \leq -\pi$. In either of those two cases, a sign flip occurs, and instead $z^{\alpha} w^{\alpha} = -(zw)^{\alpha}$.

\section{Conditions that Cause an Error in the Integrals}\label{sec:error_cond}
In (\ref{eq:f_error_combine}), a sign flip error can potentially occur in any of the three cases where variables are combined: $X^{m/2\,} \smash{\overline{Y}}^{m/2} \rightarrow (X\overline{Y})^{m/2}$, $X^{1/2\,} Y^{1/2} \rightarrow \sqrt{XY}$, or $Y^{-m/2\,} \smash{\overline{Y}}^{-m/2} \rightarrow (Y\overline{Y})^{-m/2}$. In the first and third cases, the sign flip will occur only for odd values of $m$. In the second case, $2\sqrt{XY}$ is the input to the Bessel function $I_m(z)$. However, for $\nu \in \mathbb{Z}$, $I_{\nu}(-z) = (-1)^{\nu} I_{\nu}(z)$ (which is a special case of NIST \cite[Eq.~10.34.1]{NIST_DLMF}). So, $I_m(z)$ is an even function of $z$ for even $m$ and an odd function of $z$ for odd $m$. Consequently, if a sign error occurs for $z = 2\sqrt{XY}$, it will only have an effect on the output of $I_m(z)$ if $m$ is odd. We will consider these three cases separately.

\subsection{Case 1:}\label{subsec:case1} When $X^{m/2\,} \smash{\overline{Y}}^{m/2} \neq (X\overline{Y})^{m/2}$.

For ease of reference, $X = \big[(p{+}b)+i(a{-}q)\big]/2$ and $Y = \big[(p{-}b)+i(a{+}q)\big]/2$; hence, $\overline{Y} = \big[(p{-}b)-i(a{+}q)\big]/2$. We therefore have $\theta_X = \arg X = \operatorname{atan2}(a{-}q, p{+}b)$ and $\theta_Y = \arg Y = \operatorname{atan2}(a{+}q, p{-}b)$, where $\operatorname{atan2}(y,x)$ is the ``four-quadrant'' version of the inverse tangent function $\atan(y/x)$. For $X^{m/2\,} \smash{\overline{Y}}^{m/2}$, we must generally consider $\theta_X - \theta_Y$. The exception is when $Y$ is a negative real number, where $\arg \overline{Y}=~\pi~\neq~{-}\arg Y$, which must be handled separately. We consider four cases for $\theta_X$.

\subsubsection{$\theta_X = 0$} In other words, $X$ is a positive real number. For this to occur, $a{-}q = 0$ and $p{+}b > 0$. However, in this case, it is not possible for $\theta_X + \theta_{\hspace{1pt}\overline{Y}}$ to fall outside of $(-\pi,\pi]$. Thus, no sign error will occur in this case.

\subsubsection{$\theta_X$ is undefined, i.e., $X = 0$} This occurs when $a{-}q = 0$ and $p{+}b = 0$. Trivially, though, $X^{m/2\,} \smash{\overline{Y}}^{m/2} = 0{\times}\smash{\overline{Y}}^{m/2} = 0$ and $(X\overline{Y})^{m/2} = 0^{m/2} = 0$, so there will be no sign error. (We are currently only considering odd values for $m$, so $0^0$ will not occur.) For similar reasons, a sign error will not occur if $Y=0$, so we need not consider that case.

\subsubsection{$0{\,<\,}\theta_X{\,\leq\,}\pi$} This occurs either when $a{-}q > 0$ or when $a{-}q = 0$ and $p{+}b < 0$. The only possibility for a sign error to occur in this case is when $\arg \overline{Y} > 0$. This will happen either if $a{+}q<0$, so $-\pi < \theta_Y <0$, or if $a{+}q=0$ and $p{-}b < 0$, so $Y$ is a negative real number and $\theta_Y = \theta_{\hspace{1pt}\overline{Y}} = \pi$. Considering the latter case first, if $\theta_{\hspace{1pt}\overline{Y}} = \pi$ and $\theta_X$ is strictly greater than zero, then $\theta_X + \theta_{\hspace{1pt}\overline{Y}}$ must be strictly greater than $\pi$. Hence, a condition for a sign error to occur is if either 1) $(a{-}q{\,>\,}0) \cap (a{+}q{\,=\,}0) \cap (p{-}b{\,<\,}0)$, which reduces to $(q{\,=\,}{-}a) \cap (q{\,<\,}0) \cap (p{\,<\,}b)$; or 2) $(a{-}q{\,=\,}0) \cap (p{+}b{\,<\,}0) \cap (a{+}q{\,=\,}0) \cap (p{-}b{\,<\,}0)$, which reduces to $(a{\,=\,}q{\,=\,}0) \cap (p{\,<\,}{-}|b|)$. Similarly, if $\theta_X = \pi$ (i.e., $a{-}q = 0$ and $p{+}b < 0$), then a sign error will occur for any positive value of $\theta_{\hspace{1pt}\overline{Y}}$. The case of $\theta_{\hspace{1pt}\overline{Y}} = \pi$ has already been covered above. Otherwise, a sign error will also occur when $(a{-}q{\,=\,}0) \cap (a{+}q{\,<\,}0) \cap (p{+}b{\,<\,}0)$, which reduces to $(q{\,=\,}a) \cap (q{\,<\,}0) \cap (p{\,<\,}{-}b)$.

For the remaining possibility, $\theta_X$ and $\theta_{\hspace{1pt}\overline{Y}}$ are both positive but neither equals $\pi$. This occurs when both $a{-}q > 0$ and $a{+}q < 0$, or in other words when $q < -|a|$. Under this condition, we must determine when $\theta_X + \theta_{\hspace{1pt}\overline{Y}} > \pi$. For this, we may use the following property of $\operatorname{atan2}$:\pagebreak[0]
\begin{equation}\label{eq:atan2_prop}
    \operatorname{atan2}(y,x) = 
       \begin{cases} 
       \dfrac{\pi}{2} - \atan\!\bigg(\dfrac{x}{y}\bigg), & \textrm{if } y>0;\\[10pt]
       -\dfrac{\pi}{2} - \atan\!\bigg(\dfrac{x}{y}\bigg), & \textrm{if } y<0.
       \end{cases}
\end{equation}
Then, the condition $\theta_X + \theta_{\hspace{1pt}\overline{Y}} > \pi$, or equivalently $\operatorname{atan2}(a{-}q,p{+}b) \linebreak[0]+ \operatorname{atan2}(-a{-}q,p{-}b) \linebreak[0]> \pi$, may then be rewritten as $\pi/2 - \atan\!\left(\tfrac{p{+}b}{a{-}q}\right) + \pi/2 - \atan\!\left(\tfrac{p{-}b}{-a{-}q}\right) > \pi$. After some basic algebraic manipulations, this reduces to $p < -ba/q$. It is also interesting to note that in two of the earlier cases, i.e. $(q{\,=\,}{-}a) \cap (q{\,<\,}0) \cap (p{\,<\,}b)$ and $(q{\,=\,}a) \cap (q{\,<\,}0) \cap (p{\,<\,}{-}b)$, the inequality for $p$ can be alternatively written in the same way, i.e., $p < -ba/q$ in both of those cases as well.

\subsubsection{${-}\pi{\,<\,}\theta_X{\,<\,}0$} This occurs when $a{-}q < 0$. A sign error will only occur in this case if $\arg \overline{Y}$ is also negative, so that $a{+}q > 0$. Combining the two conditions, this means $q > |a|$. The sign error will occur if $\theta_X + \theta_{\hspace{1pt}\overline{Y}} \leq {-}\pi$, or equivalently $\operatorname{atan2}(a{-}q,p{+}b) + \operatorname{atan2}(-a{-}q,p{-}b) \leq {-}\pi$. Using the lower half of (\ref{eq:atan2_prop}) and similar algebraic manipulations as before, we obtain the very similar condition $p \leq -ba/q$ (this time including equality).

To summarize, the result $X^{m/2\,} \smash{\overline{Y}}^{m/2} \neq (X\overline{Y})^{m/2}$ will occur under any of the following conditions:
\begin{subequations}\label{eq:conds_case1}
   \begin{equation}\label{eq:conds_case1a}
           \Big[(|q| > |a|) \textrm{ OR } (q = -|a| \neq 0)\Big] \textrm{ AND } p < -\frac{ba}{q};
   \end{equation}
   \begin{equation}\label{eq:conds_case1b}
           q > |a| \textrm{ AND } p = -\frac{ba}{q};
   \end{equation}
   \begin{equation}\label{eq:conds_case1c}
           q = a = 0 \textrm{ AND } p < -|b|.
   \end{equation}
\end{subequations}

\subsection{Case 2:}\label{subsec:case2}  When $X^{1/2\,} {Y}^{1/2} \neq \sqrt{XY}$.

This case is much the same as Case 1, except we are now considering $\theta_X{\,+\,}\theta_{Y}$ instead of $\theta_X{\,+\,}\theta_{\hspace{1pt}\overline{Y}}$. As such, we shall omit most of the details of the derivation. The key difference is that the inequality conditions on $a{+}q$ become the opposite of those from Case 1 (i.e., instances of ``greater than'' now become ``less than'' and vice versa). Also, with $a{+}q$ being the imaginary component of $Y$ (instead of ${-}a{-}q$ for $\overline{Y}$), the inequality for $p$ reduces to $p < -bq/a$ rather than $p < -ba/q$. Consequently, it ends up that the result $X^{m/2\,} {Y}^{m/2} \neq (XY)^{m/2}$ will occur under any of the following conditions:
\begin{subequations}\label{eq:conds_case2}
   \begin{equation}\label{eq:conds_case2a}
           \Big[(|a| > |q|) \textrm{ OR } (a = |q| \neq 0)\Big] \textrm{ AND } p < -\frac{bq}{a};
   \end{equation}
   \begin{equation}\label{eq:conds_case2b}
           a < -|q| \textrm{ AND } p = -\frac{bq}{a};
   \end{equation}
   \begin{equation}\label{eq:conds_case2c}
           a = q = 0 \textrm{ AND } p < -|b|.
   \end{equation}
\end{subequations}

\subsection{Case 3:}\label{subsec:case3} When $Y^{-m/2\,} \smash{\overline{Y}}^{-m/2} \neq (Y\overline{Y})^{-m/2}$.

Overall, it will most commonly be the case that $\theta_{\hspace{1pt}\overline{Y}} = -\theta_Y$, so $\theta_Y+\theta_{\hspace{1pt}\overline{Y}} = 0$ (or is undefined if $Y=0$). Since $0 \in (-\pi,\pi]$, there will not be a sign error. The one exception is when $Y$ is a negative real number, so $\theta_{\hspace{1pt}\overline{Y}} = \theta_Y = \pi$ and $\theta_Y+\theta_{\hspace{1pt}\overline{Y}} = 2\pi$. $Y$ will be a negative real number when $a{+}q = 0$ and $p{-}b < 0$. Hence, $Y^{-m/2\,} \smash{\overline{Y}}^{-m/2} \neq (Y\overline{Y})^{-m/2}$ will occur when:
\begin{equation}\label{eq:conds_case3}
           a = -q \textrm{ AND } p<b.
\end{equation}
Neither $Y^{-m/2}$, $\smash{\overline{Y}}^{-m/2}$, nor $(Y\overline{Y})^{-m/2}$ can be calculated when $a = -q$ and $p = b$, as $Y=0$ in that event.

\subsection{An Overall Error}\label{subsec:overall}
As seen, there are three cases that can cause a sign error when calculating the different parts of $f$. However, it could potentially occur that two of the parts produce a sign error, while the third does not. Thus, the two errors could cancel each other out, coincidentally leading to a correct overall result. Examining the conditions in (\ref{eq:conds_case1}), (\ref{eq:conds_case2}), and (\ref{eq:conds_case3}), the only time there can be an overlap is when $q = -a$, $q \leq 0$, and $a \geq 0$. If $q \neq -a$, then Case 3 will not yield a sign error, and only one of Case 1 or Case 2 can cause an error but not both. In this event, either $|q|>|a|$ (Case 1), $|a|>|q|$ (Case 2), $q = a$ and both are positive (Case 2), or $q = a$ and both are negative (Case 1). On the other hand, if $q = -a$ but $q>0$, then $a$ must be negative; likewise, if $q = -a$ but $a<0$, then $q$ must be positive. In either event, only Case 3 will yield a sign error.

When $q = -a$ and neither equals zero, then $-bq/a = -ba/q = b$. In this event, the error conditions (\ref{eq:conds_case1a}), (\ref{eq:conds_case2a}), and (\ref{eq:conds_case3}) are all equivalent to each other. Hence, either all three cases will simultaneously yield a sign error or none of them will yield an error. If all three are in error, then the overall result will also be in error.

When $a = q = 0$, then conditions (\ref{eq:conds_case1c}), (\ref{eq:conds_case2c}), and/or (\ref{eq:conds_case3}) could be satisfied. If $b \leq 0$, then the condition $p < -|b|$ is equivalent to $p<b$. Hence, again either all three cases will yield a sign error simultaneously, or none of them will yield an error. On the other hand, if $b>0$, then if $p<b$, Case 3 will yield a sign error. However, $p$ may still fall in the range of $-b \leq p < b$, which will not cause a sign error in Cases 1 and 2. Nonetheless, a single sign error will cause an overall error in the final result. If $p < -b$ with $b>0$, then Cases 1 and 2 will yield a sign error along with Case 3, and the three sign errors will again create an error in the overall result. Thus, for $a = q = 0$, it is sufficient to have $p<b$ to cause an overall error in the result.

In summary, it is not possible for two out of the three cases to yield a sign error while the third does not. Either none of the cases yields a sign error (in which event the final result is correct), a single case yields a sign error, or all three cases yield sign errors at the same time. In either of the last two events, an overall sign error will be the result.

The conditions that cause an overall sign error when calculating $f$ using (\ref{eq:f_Im_final}) can be summed up fairly succinctly. First, we define $K$ as follows:
\begin{equation}\label{eq:K_def}
    K = \begin{cases}
        q/a, & \text{if } |a| \geq |q| \textrm{ AND } a \neq 0;\\
        a/q, & \text{if } |q| \geq |a| \textrm{ AND } q \neq 0;\\
        -1, & \text{if } a = q = 0.
    \end{cases}
\end{equation}
Then, a sign error occurs under the either of the following two conditions:
\begin{subequations}\label{eq:conds_overall}
   \begin{equation}\label{eq:conds_overall_a}
           p < -bK;
   \end{equation}
   \begin{equation}\label{eq:conds_overall_b}
          \big[ (a < -|q|) \textrm{ OR } (q > |a|) \big] \textrm{ AND } p = -bK.
   \end{equation}
\end{subequations}

We lastly note that the above conditions yield a sign error specifically when calculating $f$ using (\ref{eq:f_Im_final}). However, there are a limited subset of conditions that result in either the real or imaginary part of $f$ being equal to zero while the other part is not. If (\ref{eq:conds_overall}) holds and $\mathfrak{Re}(f) = 0$, then (\ref{eq:cos_int_old}) will still correctly yield zero while (\ref{eq:sin_int_old}) will have a sign error. A straightforward way for this to happen is if $p = b = 0$; then, $B = D = 0$ and $A = -C$. With an odd value for $m$, it will then occur that $(A-iB)^{m/2}$ is purely real and $I_m\big(\sqrt{C+iD})$ is purely imaginary, or vice versa; hence, (\ref{eq:f_Im_final}) will yield a purely imaginary value. (Trivially, $Y\overline{Y}$ must be positive and real, and thus so is $(Y\overline{Y})^{-m/2}$). Likewise, if (\ref{eq:conds_overall}) holds and $\mathfrak{Im}(f) = 0$, then (\ref{eq:sin_int_old}) will still correctly yield zero while (\ref{eq:cos_int_old}) will have a sign error. A straightforward way for this latter case to occur is when $a = q = 0$; then, $B = D = 0$ and $A = C$. Consequently, $(A-iB)^{m/2}$ and $I_m\big(\sqrt{C+iD})$ will both be either purely real or purely imaginary, and so (\ref{eq:f_Im_final}) will yield a purely real value.

\section{Improved Expressions for Integrals}
It is of course possible to correct the expressions in (\ref{eq:sin_int_old}) and (\ref{eq:cos_int_old}) so that they consistently give the proper results. For example, one could multiply both expressions by $(-1)^{m}$ if the conditions in (\ref{eq:conds_overall}) hold. Alternatively, one could just not combine the alike terms in (\ref{eq:f_error_combine}). After deleting the $(\overline{Y}^{1/2})^{m\,} (\overline{Y}^{1/2})^{-m}$ terms (which cancel each other) from (\ref{eq:f_error_combine}), this would give $f = 2\pi X^{m/2\,} Y^{-m/2}  I_m\big(2\sqrt{X}\sqrt{Y}\big)$. The result of the integral in (\ref{eq:sin_int_old}) can then be expressed as $i\pi \Big[\smash{\overline{X}}^{m/2\,} \smash{\overline{Y}}^{-m/2}  I_m\Big(2\sqrt{\hspace{1pt}\overline{X}\vphantom{\bar{\bar{X}}}}\sqrt{\hspace{1pt}\overline{Y}\vphantom{\bar{\bar{Y}}}}\Big) - X^{m/2\,} Y^{-m/2}  I_m\big(2\sqrt{X}\sqrt{Y}\big)\Big]\vspace{1pt}$, and the result of the integral in (\ref{eq:cos_int_old}) can be expressed as $\pi \Big[\smash{\overline{X}}^{m/2\,} \smash{\overline{Y}}^{-m/2}  I_m\Big(2\sqrt{\hspace{1pt}\overline{X}\vphantom{\bar{\bar{X}}}}\sqrt{\hspace{1pt}\overline{Y}\vphantom{\bar{\bar{Y}}}}\Big) + X^{m/2\,} Y^{-m/2}  I_m\big(2\sqrt{X}\sqrt{Y}\big)\Big]\vspace{2pt}$.

However, we can instead derive an alternative expression that avoids the complications surrounding combining terms with non-integer exponents in the first place. Resuming the derivation by continuing from (\ref{eq:cont_int_result}):
\begin{equation}\label{eq:new_result_pt1}
\begin{aligned}
    f &= 2\pi X^m \sum_{j = 0}^{\infty} \frac{X^{j}\, Y^{j}}{j!\,(j+m)!} = 2\pi X^m \sum_{j = 0}^{\infty} \frac{X^{j}\, Y^{j}}{j!\,\Gamma(j{+}m{+}1)} \\
    &= \frac{2\pi X^m}{\Gamma(m{+}1)} \sum_{j = 0}^{\infty} \frac{X^{j\,} Y^{j\,} \Gamma(m{+}1)}{j!\,\Gamma(j{+}m{+}1)} = \frac{2\pi X^m}{m!} \sum_{j = 0}^{\infty} \frac{X^{j}\, Y^{j}}{j!\,(m{+}1)_j},
\end{aligned}
\end{equation}
where $(z)_j = \Gamma(z{+}j)/\Gamma(z)$ denotes the Pochhammer symbol. Since $j$ is a non-negative integer, $X^j$ and $Y^j$ each have a single value, and their product will equal $(XY)^{\hspace{1pt}j}$. (As a simple proof, \scalebox{0.94}{$X^j Y^j = \underbrace{X{\cdot}X{\cdot}X{\cdot}{\ldots}{\cdot}X}_{j\textrm{ copies}} \cdot \underbrace{Y{\cdot}Y{\cdot}Y{\cdot}{\ldots}{\cdot}Y}_{j\textrm{ copies}} = \underbrace{XY{\cdot}XY{\cdot}XY{\cdot}{\ldots}{\cdot}XY}_{j\textrm{ copies}}$} $= (XY)^{\hspace{1pt}j}.$) Replacing $X^{j\,} Y^j$ with $(XY)^{\hspace{1pt}j}$ in the last sum of (\ref{eq:new_result_pt1}), it can be seen that the sum is, by definition, the (confluent) hypergeometric function ${}_0F_1(;m{+}1;XY)$. This is a special case of the generalized hypergeometric function ${}_pF_q(a_1,a_2,\ldots,a_p;b_1,b_2,\ldots,b_q;z) = \sum\limits_{k=0}^{\infty} \frac{(a_1)_k (a_2)_k \cdots (a_p)_k}{(b_1)_k (b_2)_k \cdots (b_q)_k}\frac{z^k}{k!}$ (NIST \cite[Eq.~16.2.1]{NIST_DLMF}), with $p=0$ parameters in the numerator and $q = 1$ parameter in the denominator. (These are not the same ``$p$'' and ``$q$'' as otherwise used in this paper.) Therefore, we end up with:
\begin{equation}\label{eq:new_result_full}
   \begin{aligned}
    f &= \frac{2\pi X^m}{m!} {}_0F_1(;m{\,+\,}1;XY) \\
    &= \frac{2\pi}{m!} \left[ \frac{(p{+}b)+i(a{-}q)}{2} \right]^{\hspace{-1pt} m} {}_0F_1\!\left(; m{\,+\,}1; \frac{(p^2{+}q^2{-}a^2{-}b^2)+i[2(ap{\,+\,}bq)]}{4} \right).
   \end{aligned}
\end{equation}

Let us define four new constants as follows:
\begin{subequations}\label{eq:new_constants}
    \begin{equation}\label{eq:new_A}
        A^{\prime} = \frac{p+b}{2},
    \end{equation}
    \begin{equation}\label{eq:new_B}
        B^{\prime} = \frac{a-q}{2},
    \end{equation}
    \begin{equation}\label{eq:new_C}
        C^{\prime} = \frac{p^2 + q^2 - a^2 - b^2}{4},
    \end{equation}
    \begin{equation}\label{eq:new_D}
        D^{\prime} = \frac{ap + bq}{2}.
    \end{equation}
\end{subequations}
We can then express $f$ somewhat more compactly as
\begin{equation}\label{eq:new_result_compact}
    f = \frac{2\pi}{m!} (A^{\prime}+ iB^{\prime})^{m\,} {}_0F_1(; m{\,+\,}1; C^{\prime} + iD^{\prime}).
\end{equation}
We also note that, from the power series definition of ${}_0F_1\vspace{2pt}$ and the property $\overline{z^{\alpha}} = \overline{z}^{\hspace{1pt}\alpha}$ for $\alpha \in \mathbb{R}$, it follows that ${}_0F_1(;b_1;\overline{z}) = \overline{{}_0F_1(;b_1;z)}$ for $b_1 \in \mathbb{R}$. Using $(f{-}\overline{f})/(2i)$ and $(f{+}\overline{f})/2$, we ultimately arrive at the new expressions for the integrals:
\begin{align}\label{eq:sin_int_new}
\hspace{24pt}&\int_{0}^{2\pi} \exp(p \cos x + q \sin x) \sin(a \cos x + b \sin x - mx) \, dx \nonumber \\[-10pt]
&\quad= \frac{i\pi}{m!} \big[(A^{\prime}{\hspace{1pt}-\,}iB^{\prime})^{m\,} {}_0F_1(; m{\,+\,}1; C^{\prime}{\hspace{1pt}-\,}iD^{\prime}) - (A^{\prime}{\hspace{1pt}+\,}iB^{\prime})^{m\,} {}_0F_1(; m{\,+\,}1; C^{\prime}{\hspace{1pt}+\,}iD^{\prime})\big]
\end{align}
\begin{align}\label{eq:cos_int_new}
\hspace{24pt}&\int_{0}^{2\pi} \exp(p \cos x + q \sin x) \cos(a \cos x + b \sin x - mx) \, dx \nonumber \\[-10pt]
&\quad= \frac{\pi}{m!} \big[(A^{\prime}{\hspace{1pt}-\,}iB^{\prime})^{m\,} {}_0F_1(; m{\,+\,}1; C^{\prime}{\hspace{1pt}-\,}iD^{\prime}) + (A^{\prime}{\hspace{1pt}+\,}iB^{\prime})^{m\,} {}_0F_1(; m{\,+\,}1; C^{\prime}{\hspace{1pt}+\,}iD^{\prime})\big]
\end{align}
There is one corner case with the new expressions that also existed in the originals. In the event $A^{\prime} = B^{\prime} = m = 0$, then the expressions contain $0^0$. In this event, the terms $(A^{\prime} + iB^{\prime})^{m}$ and $(A^{\prime} - iB^{\prime})^{m}$ should be treated as $\lim_{z \rightarrow 0} z^0 = 1$.

The expressions in (\ref{eq:sin_int_new}) and (\ref{eq:cos_int_new}) offer several improvements over the original expressions in (\ref{eq:sin_int_old}) and (\ref{eq:cos_int_old}). Most importantly, they avoid the sign errors that can occur with the original expressions. The new expressions are also somewhat more compact, having fewer terms to calculate than the original ones. This is mostly since the term $\big[(b{-}p)^2 + (a{+}q)^2\big]^{-m/2}$ is no longer present. Because of the absence of that term, (\ref{eq:sin_int_new}) and (\ref{eq:cos_int_new}) can therefore also be used in the event $(b{-}p)^2 + (a{+}q)^2 = 0$, which is a limitation of the original expressions.

There is one further advantage to the newer expressions. So far in this paper, it has been implicitly assumed that $a$, $b$, $p$, and $q$ are all real-valued constant parameters. However, there is nothing in particular in the derivations of (\ref{eq:orig_int})--(\ref{eq:f_manip}) and (\ref{eq:new_result_pt1})--(\ref{eq:new_result_compact}) that prevents those constants from being complex numbers instead. Thus, (\ref{eq:new_result_compact}) also will give the correct result for the integral in (\ref{eq:orig_int}) when using complex-valued constants. However, there is now the complication that the solution to the integral in (\ref{eq:sin_int_new}) would no longer be obtained from $(f{-}\overline{f})/(2i)$, and the solution to the integral in (\ref{eq:cos_int_new}) would no longer be obtained from $(f{+}\overline{f})/2$. Instead, those two integrals must be reworked a bit first. For brevity of notation, we shall denote the real part of a constant ``$u$'' (i.e., $\mathfrak{Re}\,u$) by $u_R$ and the imaginary part (i.e., $\mathfrak{Im}\,u$) by $u_I$. For (\ref{eq:sin_int_new}), we have:
\begin{equation}\label{eq:sin_int_cplx_init}
\begin{aligned}
&\int_{0}^{2\pi} \exp\!\big[(p_R{+}i\hspace{1pt}p_I) \cos x + (q_R{+}i\hspace{1pt}q_I) \sin x\big]  \\
&\qquad \times \sin\!\big[(a_R{+}i\hspace{1pt}a_I) \cos x + (b_R{+}i\hspace{1pt}b_I) \sin x - mx\big] \, dx  \\
&=\frac{1}{2i}\int_{0}^{2\pi} \exp\!\big[(p_R{+}i\hspace{1pt}p_I) \cos x + (q_R{+}i\hspace{1pt}q_I) \sin x\big]  \\ 
& \hspace{60pt}\times \Big\{\exp\!\Big(i\big[(a_R{+}i\hspace{1pt}a_I) \cos x + (b_R{+}i\hspace{1pt}b_I) \sin x - mx\big]\Big)  \\
& \hspace{86pt} {}- \exp\!\Big({-}i\big[(a_R{+}i\hspace{1pt}a_I) \cos x + (b_R{+}i\hspace{1pt}b_I) \sin x - mx\big]\Big) \Big\} \, dx  \\
&=\frac{i}{2}\times\underbrace{
    \begin{matrix}\hspace{-52pt}\displaystyle{\int_{0}^{2\pi}} \exp\!\Big(\big[(p_R{+}a_I)+i\hspace{1pt}(p_I{-}a_R)\big] \cos x \\
       \hspace{52pt}+ \big[(q_R{+}b_I)+i\hspace{1pt}(q_I{-}b_R)\big] \sin x\Big)e^{+imx} \, dx
    \end{matrix}}_{f_1}  \\ 
&\quad - \frac{i}{2}\times\underbrace{
        \begin{matrix}\hspace{-52pt}\displaystyle{\int_{0}^{2\pi}} \exp\!\Big(\big[(p_R{-}a_I)+i\hspace{1pt}(p_I{+}a_R)\big] \cos x \\
        \hspace{52pt}+ \big[(q_R{-}b_I)+i\hspace{1pt}(q_I{+}b_R)\big] \sin x\Big)e^{-imx} \, dx
        \end{matrix}}_{f_2}  \\ 
& = \frac{i}{2} (f_1 - f_2).
\end{aligned}
\end{equation}
For (\ref{eq:cos_int_new}), we end up with $(f_1 + f_2)/2$. One can observe that if $p$, $q$, $a$, and $b$ are all real values, so $p_I = q_I = a_I = b_I = 0$, $f_2$ simplifies to $f$ in (\ref{eq:orig_int}), and $f_1$ reduces to $\overline{f}$.

Hence, with $f_1$ and $f_2$, we have two expressions that are extremely similar to that of (\ref{eq:orig_int}). Thus, the results of the integrals are still given in the form of (\ref{eq:new_result_compact}). The key difference is simply that original constants $p$, $q$, $a$, and $b$ are replaced by the real part of that constant plus or minus the imaginary part of a different constant; for example, $p$ becomes $p_R{+}a_I$ or $p_R{-}a_I$. (The two ``$\cos$'' constants $p$ and $a$ and the two ``$\sin$'' constants $q$ and $b$ end up paired together.) $f_1$ also has the slight complication that it contains a ``${+}m$'' instead of a ``${-}m$''. However, we can instead consider $\overline{f_1}$ to convert that plus to a minus, then undo the conjugate for the final result. It is then simply a matter of making the appropriate substitution of constants in (\ref{eq:new_result_compact}) to obtain the result. The final expressions for the integrals are as follows:\pagebreak[0]
\begin{align}\label{eq:sin_int_cplx} 
\hspace{24pt}&\int_{0}^{2\pi} \!\!\exp(p \cos x{\,+\,}q \sin x) \sin(a \cos x{\,+\,}b \sin x{\,-\,}mx)\hspace{1pt}dx \hspace{5pt}\scalebox{0.875}{$[p, q, a, \textrm{and } b \textrm{ complex-valued}]$}\nonumber \\[-10pt]
&= \frac{i\pi}{m!} \big[(A_1{\hspace{1pt}+\,}iB_1)^{m\,} {}_0F_1(; m{\,+\,}1; C_1{\hspace{1pt}+\,}iD_1) - (A_2{\hspace{1pt}+\,}iB_2)^{m\,} {}_0F_1(; m{\,+\,}1; C_2{\hspace{1pt}+\,}iD_2)\big]
\end{align}\pagebreak[0]
\begin{align}\label{eq:cos_int_cplx}
\hspace{24pt}&\int_{0}^{2\pi} \!\!\exp(p \cos x{\,+\,}q \sin x) \cos(a \cos x{\,+\,}b \sin x{\,-\,}mx)\hspace{1pt}dx \hspace{5pt}\scalebox{0.875}{$[p, q, a, \textrm{and } b \textrm{ complex-valued}]$} \nonumber \\[-10pt]
&= \frac{\pi}{m!} \big[(A_1{\hspace{1pt}+\,}iB_1)^{m\,} {}_0F_1(; m{\,+\,}1; C_1{\hspace{1pt}+\,}iD_1) + (A_2{\hspace{1pt}+\,}iB_2)^{m\,} {}_0F_1(; m{\,+\,}1; C_2{\hspace{1pt}+\,}iD_2)\big]
\end{align}
where
\begin{subequations}\label{eq:cplx_params}
\begin{equation}\label{eq:cplx_params_A1}
        A_1 = \dfrac{p_R + a_I - q_I + b_R}{2},
\end{equation}
\begin{equation}\label{eq:cplx_params_A2}
        A_2 = \dfrac{p_R - a_I + q_I + b_R}{2},
\end{equation}
\begin{equation}\label{eq:cplx_params_B1}
        B_1 = \dfrac{p_I - a_R + q_R  + b_I}{2},
\end{equation}
\begin{equation}\label{eq:cplx_params_B2}
        B_2 = \dfrac{p_I + a_R - q_R + b_I}{2},
\end{equation}
\begin{equation}\label{eq:cplx_params_C1}
        C_1 = \dfrac{(p_R{\,+\,}a_I)^2 + (q_R{\,+\,}b_I)^2 - (p_I{\,-\,}a_R)^2 - (q_I{\,-\,}b_R)^2}{4},
\end{equation}
\begin{equation}\label{eq:cplx_params_C2}
        C_2 = \dfrac{(p_R{\,-\,}a_I)^2 + (q_R{\,-\,}b_I)^2 - (p_I{\,+\,}a_R)^2 - (q_I{\,+\,}b_R)^2}{4},
\end{equation}
\begin{equation}\label{eq:cplx_params_D1}
        D_1 = \dfrac{(p_I{\,-\,}a_R)(p_R{\,+\,}a_I) + (q_I{\,-\,}b_R)(q_R{\,+\,}b_I)}{2},
\end{equation}
\begin{equation}\label{eq:cplx_params_D2}
        D_2 = \dfrac{(p_I{\,+\,}a_R)(p_R{\,-\,}a_I) + (q_I{\,+\,}b_R)(q_R{\,-\,}b_I)}{2}.
\end{equation}
\end{subequations}
Again, if $p_I = q_I = a_I = b_I = 0$ so that the constants are real-valued, then the expressions in (\ref{eq:cplx_params}) will reduce so that $A_1 = A_2 = A^{\prime}$, $B_1 = {-}B^{\prime}$, $B_2 = B^{\prime}$, $C_1 = C_2 = C^{\prime}$, $D_1 = {-}D^{\prime}$, and $D_2 = D^{\prime}$. Hence, (\ref{eq:sin_int_cplx}) and (\ref{eq:cos_int_cplx}) will reduce to (\ref{eq:sin_int_new}) and (\ref{eq:cos_int_new}), respectively.

We lastly note that complex-valued constants can also be technically used with the form of $f$ given by (\ref{eq:f_Im_final}), with its restriction changing from $(b-p)^2 + (a+q)^2 > 0$ to $(b-p)^2 + (a+q)^2 \neq 0$. However, the same type of sign errors will still result from its use, and the analysis of the conditions where those sign errors occur would be considerably more complicated. Use of (\ref{eq:new_result_compact}) therefore still remains the better option.

\section{Comments on Other Integrals}\label{sec:other_comments}
The results of (\ref{eq:sin_int_new}) and (\ref{eq:cos_int_new}) can be used to verify the results of a few other related integrals in Gradshteyn and Ryzhik \cite{gr14}. Specifically, we consider the following integrals, which are special cases of (\ref{eq:sin_int_new}) and (\ref{eq:cos_int_new}). (In some instances, we add accents below to aid in distinguishing between constants.) 

\subsection{}\label{subsec:3.931_4} \textbf{Integral 1} (Gradshteyn and Ryzhik \cite[Eq.~3.931~4]{gr14}):
\begin{equation}\label{eq:3.931_4}
\int_{0}^{\pi} e^{-p^{\prime} \cos x} \cos(p^{\prime} \sin x) \, dx = \frac{1}{2}\int_{0}^{2\pi} e^{-p^{\prime} \cos x} \cos(p^{\prime} \sin x) \, dx = \pi
\end{equation}

The equality between the two integrals can be determined by noting the value of $\sin x$ from $\pi$ to $2\pi$ equals the negative of the value from $\pi$ to 0, whereas the value of $\cos x$ from $\pi$ to $2\pi$ equals the same value from $\pi$ to 0. Thus, the equation takes on the same values between $\pi$ to $2\pi$ as it does from 0 to $\pi$, albeit in the reverse order. Hence, the value of the integral over those two ranges of angles is the same.

The integral on the right corresponds to (\ref{eq:cos_int_new}) with $p = -p^{\prime}$, $b = p^{\prime}$, and $q = a = m = 0$. We therefore have $A^{\prime}$, $B^{\prime}$, $C^{\prime}$, and $D^{\prime}$ all equal to zero. Consequently, as stated before, $(A^{\prime} \pm iB^{\prime})^{m}$ should be treated as being equal to 1. (\ref{eq:cos_int_new}) then gives $(\pi/0!)\big[1{\cdot}{\,}_0F_1(;1;0)+1{\cdot}{\,}_0F_1(;1;0)\big] = \pi[1{+}1] = 2\pi$. Therefore, after multiplying by the leading factor of $1/2$, this matches with and confirms the original result. This result also applies in the degenerate case of $p^{\prime} = 0$: $\int_{0}^{\pi} e^0 \cos(0)\,dx = \int_{0}^{\pi} 1{\,\cdot\,}1\,dx = \pi$. We furthermore note that the same result is achieved in the event $p^{\prime}$ is complex. In this case, all the constants in (\ref{eq:cplx_params}) work out to be 0, so (\ref{eq:cos_int_cplx}) yields the value as (\ref{eq:cos_int_new}). It is also worth noting that the same result holds if the $p^{\prime}$ inside the $\cos$ term is replaced by $-p^{\prime}$, since $\cos(-z) = \cos(z)$. By extension, the minus sign in the $\exp$ term may also be removed (e.g., by setting $p^{\prime} = -u$). Hence, \cite[Eq.~3.931~4]{gr14} can be generalized to be $\int_{0}^{\pi} e^{\pm p^{\prime} \cos x} \cos(p^{\prime} \sin x) \, dx = (1/2)\int_{0}^{2\pi} e^{\pm p^{\prime} \cos x} \cos(p^{\prime} \sin x) \, dx = \pi$.

\subsection{}\label{subsec:3.932_1_2} \textbf{Integrals 2 and 3} (Gradshteyn and Ryzhik \cite[Eqs.~3.932~1 and 3.931~2]{gr14}):
\begin{equation}\label{eq:3.932_1}
\int_{0}^{\pi} \!e^{p^{\prime} \cos x} \sin(p^{\prime} \sin x) \sin mx \, dx = \frac{1}{2} \int_{0}^{2\pi} \!e^{p^{\prime} \cos x} \sin(p^{\prime} \sin x) \sin mx \, dx = \frac{\pi p^{\prime \, m}}{2\,m!}\hspace{-2pt}
\end{equation}
\begin{equation}\label{eq:3.932_2}
\int_{0}^{\pi} \!e^{p^{\prime} \cos x} \cos(p^{\prime} \sin x) \cos mx \, dx = \frac{1}{2} \int_{0}^{2\pi} \!e^{p^{\prime} \cos x} \cos(p^{\prime} \sin x) \cos mx \, dx = \frac{\pi p^{\prime \, m}}{2\,m!}\hspace{-6pt}
\end{equation}

The equality between the two integrals again is a result of the mirror symmetry of the equation between 0 to $\pi$ and between $\pi$ to $2\pi$; hence, the integral over those two ranges of angles gives the same value. The expressions in both equations are not quite in the required form. However, we can make use of the trigonometric product identities $\sin \theta \cdot \sin \phi = \big[\!\cos(\theta{-}\phi)-\cos(\theta{+}\phi)\big]/2$ and $\cos \theta \cdot \cos \phi = \big[\!\cos(\theta{-}\phi)+\cos(\theta{+}\phi)\big]/2$. The expression in (\ref{eq:3.932_1}) can then be rewritten as
\begin{equation}\label{eq:3.932_1_alt}
\frac{1}{4} \int_{0}^{2\pi} \big[e^{p^{\prime} \cos x} \cos(p^{\prime} \sin x - mx) -  e^{p^{\prime} \cos x} \cos(p^{\prime} \sin x + mx)\big] \, dx,
\end{equation}
and the expression in (\ref{eq:3.932_2}) can be rewritten as
\begin{equation}\label{eq:3.932_2_alt}
\frac{1}{4} \int_{0}^{2\pi}  \big[e^{p^{\prime} \cos x} \cos(p^{\prime} \sin x - mx) +  e^{p^{\prime} \cos x} \cos(p^{\prime} \sin x + mx)\big] \, dx.
\end{equation}
The last $\cos$ term in both (\ref{eq:3.932_1_alt}) and (\ref{eq:3.932_2_alt}) can be written equivalently as $\cos(-p^{\prime} \sin x - mx)$.

We therefore end up with two applications of (\ref{eq:cos_int_new}); in the first, $p = p^{\prime}$, $b=p^{\prime}$, and $a = q = 0$, while in the second, $b = -p^{\prime}$, and the other constants are the same as in the first. This yields $A^{\prime} = p^{\prime}$ in the first case and $A^{\prime} = 0$ in the second case. In both cases, $B^{\prime} = C^{\prime} = D^{\prime} = 0$. Substituting these values into (\ref{eq:cos_int_new}) yields $2\pi p^{\prime\,m}/m!$ for the first half of the integrals in (\ref{eq:3.932_1_alt}) and (\ref{eq:3.932_2_alt}). The second half, however, depends on whether $m>0$ or $m=0$, on account of $0^m$ terms. If $m>0$, then $0^m = 0$, and the second half of the integral reduces to zero. On the other hand, if $m=0$, then $0^0$ should be treated as equal to 1, and the result for the second half is $2\pi$ (the same as seen in the previous section). Consequently, the result for (\ref{eq:3.932_1}) and (\ref{eq:3.932_1_alt}) is $\big[2\pi p^{\prime\,m}/m! - 0\big]/4 = (\pi p^{\prime\,m})/(2\,m!)$ when $m>0$, but $\big[2\pi p^{\prime\,m}/m! - 2\pi\big]/4 = 0$ when $m=0$. Likewise, the the result for (\ref{eq:3.932_2}) and (\ref{eq:3.932_2_alt}) is $\big[2\pi p^{\prime\,m}/m! + 0\big]/4 = (\pi p^{\prime\,m})/(2\,m!)$ when $m>0$, but $\big[2\pi p^{\prime\,m}/m! + 2\pi\big]/4 = \pi$ when $m=0$. Hence, it must be specified that the original results for (\ref{eq:3.932_1}) and (\ref{eq:3.932_2}) are \underline{only applicable when $m>0$}. The results for $m=0$ can be cross-checked by substituting $m=0$ into the original integrals. For (\ref{eq:3.932_1}), $\sin mx = \sin 0 = 0$, and thus $\int_{0}^{\pi} 0 \,dx = 0$. For (\ref{eq:3.932_2}), $\cos mx = \cos 0 = 1$, and thus the integral reduces to the equation in Section \ref{subsec:3.931_4}. As already seen, the result of that integral has been confirmed to be $\pi$.

In considering the case of complex-valued $p^{\prime}$, for the first half of the integrals in (\ref{eq:3.932_1_alt}) and (\ref{eq:3.932_2_alt}), from (\ref{eq:cplx_params}) we obtain $A_1 = A_2 = p_R^{\prime}$, $B_1 = B_2 = p_I^{\prime}$, and the remaining constants equal 0. Thus, $A_1{\,+\,}iB_1 = A_2{\,+\,}iB_2 = p_R^{\prime}{\,+\,}ip_I^{\prime}$, which is simply $p^{\prime}$. The first half of the integral thus evaluates to $2\pi p^{\prime\,m}/m!$, the same as if $p^{\prime}$ is real-valued. Likewise, for the second half of (\ref{eq:3.932_1_alt}) and (\ref{eq:3.932_2_alt}), all the constants in (\ref{eq:cplx_params}) are equal to zero. We therefore get the same results as the real case: 0 if $m>0$, and $2\pi$ if $m=0$. Hence, the overall solutions to (\ref{eq:3.932_1})/(\ref{eq:3.932_1_alt}) and to (\ref{eq:3.932_2})/(\ref{eq:3.932_2_alt}) remain the same whether $p^{\prime}$ is real-valued or complex-valued.

\subsection{}\label{subsec:3.936_1} \textbf{Integral 4} (Gradshteyn and Ryzhik \cite[Eq.~3.936~1]{gr14}):
\begin{equation}\label{eq:3.936_1}
\int_{0}^{2\pi}  e^{p^{\prime} \cos x} \cos(p^{\prime} \sin x - mx) \, dx = 2\int_{0}^{\pi}  e^{p^{\prime} \cos x} \cos(p^{\prime} \sin x - mx) \, dx = \frac{2\pi p^{\prime \, m}}{m!}
\end{equation}

This integral is the same as we considered in the previous section for the first half of (\ref{eq:3.932_1_alt}) and (\ref{eq:3.932_2_alt}). We have confirmed that $2\pi p^{\prime\,m}/m!$ is correct for all non-negative integers $m$, whether $p$ is real or complex.

\subsection{}\label{subsec:3.936_2_3} \textbf{Integrals 5 and 6} (Gradshteyn and Ryzhik \cite[Eqs.~3.936~2 and 3.936~3]{gr14}):
\begin{equation}\label{eq:3.936_2}
\int_{0}^{2\pi}  e^{p^{\prime} \sin x} \sin(p^{\prime} \cos x + mx) \, dx  = \frac{2\pi p^{\prime \, m}}{m!} \sin\!\left(\frac{m \pi}{2}\right) \qquad [p^{\prime}>0]
\end{equation}
\begin{equation}\label{eq:3.936_3}
\int_{0}^{2\pi}  e^{p^{\prime} \sin x} \cos(p^{\prime} \cos x + mx) \, dx  = \frac{2\pi p^{\prime \, m}}{m!} \cos\!\left(\frac{m \pi}{2}\right) \qquad [p^{\prime}>0]
\end{equation}

To begin, we rewrite the formulas as $(-1)e^{p^{\prime} \sin x} \sin(-p^{\prime} \cos x - mx)$ and $e^{p^{\prime} \sin x} \times \cos(-p^{\prime} \cos x - mx)$ to obtain the required ${-}mx$ instead of ${+}mx$. For both formulas, we then have $q = p^{\prime}$, $a = {-}p^{\prime}$, and $p = b = 0$. These parameters give $B^{\prime} = -p^{\prime}$ and $A^{\prime} = C^{\prime} = D^{\prime} = 0$. Substituting these values into (\ref{eq:sin_int_new}) and multiplying by the $-1$ factor in front of the rewritten formula yields
\begin{equation}\label{eq:3.936_2_pt1}
\begin{aligned}    
    &(-1) \frac{i\pi}{m!} \big[(0+ip^{\prime})^m {}_0F_1(;m{+}1;0) - (0 - ip^{\prime})^m {}_0F_1(;m{+}1;0)\big] \\
    &\quad = \frac{\pi}{i\,m!} \big[(ip^{\prime})^m - (-ip^{\prime})^m \big].
\end{aligned}
\end{equation}

Since $m$ is a non-negative integer, the factor of $p^{\prime \, m}$ may safely be pulled out of both terms. We also can express $\pm i$ as $e^{\pm i \pi/2}$. This gives
\begin{equation}\label{eq:3.936_2_confirm}
    \frac{\pi p^{\prime \, m}}{i\,m!} \big[(i)^m - (-i)^m \big] = \frac{2 \pi p^{\prime \, m}}{m!} \frac{e^{im{\pi}/2} - e^{-im{\pi}/2}}{2i} = \frac{2\pi p^{\prime \, m}}{m!} \sin\!\left(\frac{m \pi}{2}\right).
\end{equation}
Similarly, substituting the constants into (\ref{eq:cos_int_new}) and following the same steps ultimately gives
\begin{equation}\label{eq:3.936_3_confirm}
    \frac{2 \pi p^{\prime \, m}}{m!} \frac{e^{im{\pi}/2} + e^{-im{\pi}/2}}{2} = \frac{2\pi p^{\prime \, m}}{m!} \cos\!\left(\frac{m \pi}{2}\right).
\end{equation}
This therefore leads to a significant observation: the results of (\ref{eq:3.936_2_confirm}) and (\ref{eq:3.936_3_confirm}) are valid even for $p^{\prime}<0$ and for $p^{\prime}=0$ if $m>0$. They also hold for $p^{\prime} = m = 0$ if $p^{\prime \, m}$ is treated as equal to one. Consequently, in the original results given in \cite[Eqs.~3.936~2 and 3.936~3]{gr14}, \underline{the restriction $p^{\prime}>0$ is unnecessary}.

Extending the consideration to complex-valued $p^{\prime}$, substituting the parameters into (\ref{eq:cplx_params}) gives $A_1 = -p_{I}^{\prime}$, $A_2 = p_{I}^{\prime}$, $B_1 = p_{R}^{\prime}$, $B_2 = -p_{R}^{\prime}$, and $C_1 = C_2 = D_1 = D_2 = 0$. Inserting these values into (\ref{eq:sin_int_cplx}) and again multiplying by the $-1$ factor yields
\begin{equation}\label{eq:3.936_2_complex}
\begin{aligned}
    &(-1) \frac{i\pi}{m!} \big[(-p_{I}^{\prime}+ip_{R}^{\prime})^m  - (p_{I}^{\prime} - ip_{R}^{\prime})^m \big] \\
    &\quad= \frac{\pi}{i\,m!} \Big(\big[i(p_{R}^{\prime} + ip_{I}^{\prime})\big]^m - \big[(-i)(p_{R}^{\prime} + ip_{I}^{\prime})\big]^m \Big) \\
    &\quad= \frac{\pi}{i\,m!} \big[(ip^{\prime})^m - (-ip^{\prime})^m \big].
\end{aligned}
\end{equation}
This is the same as the end of (\ref{eq:3.936_2_pt1}). Although $p^{\prime}$ is now complex-valued, we can still safely pull out the factor of $p^{\prime \, m}$ as before. Ultimately, we find the same result as in (\ref{eq:3.936_2_confirm}) holds for complex-valued $p^{\prime}$, as does the result in (\ref{eq:3.936_3_confirm}).

\subsection{}\label{subsec:3.936_4} \textbf{Integral 7} (Gradshteyn and Ryzhik \cite[Eq.~3.936~4]{gr14}):
\begin{equation}\label{eq:3.936_4}
\int_{0}^{2\pi}  e^{\cos x} \sin(mx - \sin x) \, dx = 0
\end{equation}

After rewriting the formula as $(-1)e^{\cos x} \sin(\sin x - mx)$, we find this is a special case of (\ref{eq:sin_int_new}) with $p = b = 1$ and $q = a = 0$. This therefore gives $A^{\prime} = 1$ and $B^{\prime} = C^{\prime} = D^{\prime} = 0$. Substituting these into (\ref{eq:sin_int_new}) yields the first half of the formula equal to the second half. Subtracting them hence gives zero, confirming the original result.

In a sense, this may be considered to be a special case of the ``$\sin$'' counterpart to the ``$\cos$'' integral in the first half of (\ref{eq:3.932_1_alt}) and (\ref{eq:3.932_2_alt}), in which $p^{\prime} = 1$. In fact, the same result would be obtained in any case where $b = p$ and $q = a = 0$, including complex values. Thus, \cite[Eq.~3.936~4]{gr14} could be generalized to $\int_{0}^{2\pi}  e^{p \cos x} \sin(mx - p\sin x) \, dx = 0$ or $\int_{0}^{2\pi}  e^{p \cos x} \sin(p\sin x - mx) \, dx = 0$. Furthermore, we may consider $e^{p\cos x} \sin(p\sin x + mx)$ instead of $e^{p\cos x} \sin(p\sin x - mx)$. In this case, we end up with $b = -p$ and $q = a = 0$. As seen already in earlier subsections, this makes all the constants in (\ref{eq:new_constants}) and (\ref{eq:cplx_params}) equal to zero. Consequently, the first and second parts of (\ref{eq:sin_int_new}) equal each other, as do the first and second parts of (\ref{eq:sin_int_cplx}). Subtracting the two parts thus again equals zero. Hence, \cite[Eq.~3.936~4]{gr14} could be generalized even further to $\int_{0}^{2\pi}  e^{p \cos x} \sin(p\sin x \pm mx) \, dx = 0$, which applies for both real-valued and complex-valued $p$.

\subsection{}\label{subsec:3.937_3_4} \textbf{Integrals 8 and 9} (Gradshteyn and Ryzhik \cite[Eqs.~3.937~3 and 3.937~4]{gr14}):
\begin{equation}\label{eq:3.937_3}
\begin{aligned}
\int_{0}^{2\pi}  &\!\exp(p \cos x + q \sin x) \sin(q \cos x - p \sin x + mx) \, dx \\
&= \frac{2\pi}{m!} \big( p^{2} + q^{2}\big)^{m/2} \sin\!\left(\! m \atan \frac{q}{p}\right)
\end{aligned}
\end{equation}
\begin{equation}\label{eq:3.937_4}
\begin{aligned}
\int_{0}^{2\pi}  &\!\exp(p \cos x + q \sin x) \cos(q \cos x - p \sin x + mx) \, dx  \\
&= \frac{2\pi}{m!} \big( p^{2} + q^{2}\big)^{m/2} \cos\!\left(\! m \atan \frac{q}{p}\right)
\end{aligned}
\end{equation}

Once again we begin by rewriting the formulas as $(-1)\exp(p \cos x + q \sin x) \times \sin(-q \cos x + p \sin x  - mx)$ and $\exp(p \cos x + q \sin x) \cos(-q \cos x + p \sin x - mx)$. We thus have a special case of (\ref{eq:sin_int_new}) and (\ref{eq:cos_int_new}) where $a = -q$ and $b = p$. These parameters give $A^{\prime} = p$, $B^{\prime} = -q$, and $C^{\prime} = D^{\prime} = 0$. Substituting these constants into (\ref{eq:sin_int_new}), multiplying by the factor of $-1$ in front, then doing a bit of algebraic manipulation gives  $\big[\pi/(i\,m!)\big]{\cdot}\big[(p{\,+\,}iq)^m-(p{\,-\,}iq)^m\big]$. Let $z = p+iq$, so that the expression may be written as
\begin{equation}\label{eq:3.937_3_deriv}
    \frac{\pi}{i\,m!} \big[z^m-\overline{z}^{\hspace{1pt}m}\big] = \frac{2\pi}{m!}\, |z|^m \,\frac{e^{i\,m\arg z}-e^{-i\,m\arg z}}{2i} = \frac{2\pi}{m!} \big(p^2{\,+\,}q^2\big)^{m/2} \sin(m \arg z)
\end{equation}

The result in (\ref{eq:3.937_3_deriv}) is similar to the original result in (\ref{eq:3.937_3}). However, it must be noted that $\arg(p+iq)$ does not equal $\atan(q/p)$ in all circumstances. The former yields angles in the range of $(-\pi,\pi]$ radians, whereas the latter only yields angles in $(-\pi/2,\pi/2)$. If $p < 0$, then there will be a difference of $\pm\pi$ between the two. In general, one may say that $\arg(p+iq) = \atan(q/p) + k\pi$, where $k$ may equal $0$, $1$, or $-1$ depending on the values of $p$ and $q$. We may substitute this into (\ref{eq:3.937_3_deriv}), then make make use of the trigonometric identity $\sin(\theta+\phi) = \sin \theta \cos \phi + \cos \theta \sin \phi$:\pagebreak[0]
\begin{equation}
\begin{aligned}
    &\frac{2\pi}{m!} \big(p^2{\,+\,}q^2\big)^{m/2} \sin\!\big[m \arg (p+iq)\big] = \frac{2\pi}{m!} \big(p^2{\,+\,}q^2\big)^{m/2} \sin\!\left[m \!\left(\!\atan \frac{q}{p} + k\pi\!\right)\right]  \\
    &\quad = \frac{2\pi}{m!} \big(p^2{\,+\,}q^2\big)^{m/2} \left[\sin\!\left(\!m \atan \frac{q}{p}\right)\cos(km\pi) + \cos\!\left(\!m \atan \frac{q}{p}\right)\sin(km\pi)\right]
\end{aligned}
\end{equation}
$\sin(km\pi)$ will equal 0 for integer values of $k$ and $m$; thus, the second part of the above equation will vanish. On the other hand, for $\cos(km\pi)$, if either $m$ is even or $k = 0$, then we are taking the cosine of an integer multiple of $2\pi$, which yields a value of 1. However, if $m$ is odd and $k = \pm 1$, then we are taking the cosine of an odd integer multiple of $\pi$, which yields a value of $-1$. Hence, the original formula in \cite[Eq.~3.937~3]{gr14} will \underline{yield a sign error if $m$ is odd and $p<0$}. (The formula also might not be able to be used depending on whether the method of calculating $\atan(\cdot)$ accepts a number divided by zero as an input.)

Similarly, substituting the constants into (\ref{eq:cos_int_new}) will yield $(2\pi/m!) \big(p^2{\,+\,}q^2\big)^{m/2} \times \cos(m \arg z)$. We may then use the trigonometric identity $\cos(\theta + \phi) = \cos \theta \cos \phi - \sin \theta \sin \phi$ to convert $\cos(m \arg z)$ to $\cos\!\big[m \atan(q/p)\big]{\cdot}\cos(km\pi) - \sin\!\big[m \atan(q/p)\big] \times \sin(km\pi)$. Consequently, it can be seen that the original formula in \cite[Eq.~3.937~4]{gr14} will also \underline{yield a sign error if $m$ is odd and $p<0$}.

Both formulas may be corrected by multiplying by $(-1)^m$ if $p<0$. However, one can alternatively express \cite[Eq.~3.937~3]{gr14} in either of the following ways, which both corrects the error and allows for the use of $p=0$:
\begin{subequations}\label{eq:3.937_3_correct}
\begin{align}
\int_{0}^{2\pi}  &\exp(p \cos x + q \sin x) \sin(q \cos x - p \sin x + mx) \, dx  \nonumber \\
&\qquad = \frac{2\pi}{m!} \big( p^{2} + q^{2}\big)^{m/2} \sin\!\big[m \operatorname{atan2}(q,p)\big] \label{eq:3.937_3_correct_v1}\\
&\qquad = \frac{2\pi}{m!} \big( p^{2} + q^{2}\big)^{m/2} \sin\!\big[m \arg(p+iq)\big] .\label{eq:3.937_3_correct_v2}
\end{align}
\end{subequations}
Gradshteyn and Ryzhik \cite[Eq.~3.937~4]{gr14} may be expressed correctly in a similar way, replacing $\sin$ in (\ref{eq:3.937_3_correct_v1}) and (\ref{eq:3.937_3_correct_v2}) with $\cos$:
\begin{subequations}\label{eq:3.937_4_correct}
\begin{align}
\int_{0}^{2\pi}  &\exp(p \cos x + q \sin x) \cos(q \cos x - p \sin x + mx) \, dx  \nonumber \\
&\qquad = \frac{2\pi}{m!} \big( p^{2} + q^{2}\big)^{m/2} \cos\!\big[m \operatorname{atan2}(q,p)\big] \label{eq:3.937_4_correct_v1}\\
&\qquad = \frac{2\pi}{m!} \big( p^{2} + q^{2}\big)^{m/2} \cos\!\big[m \arg(p+iq)\big] .\label{eq:3.937_4_correct_v2}
\end{align}
\end{subequations}

When considering complex-valued $p$ and $q$, substitution of the parameters into (\ref{eq:cplx_params}) gives $A_1 = p_R-q_I$, $A_2 = p_R+q_I$, $B_1 = p_I+q_R$, $B_2 = p_I-q_R$, and $C_1 = C_2 = D_1 = D_2 = 0$. Inserting these constants into (\ref{eq:sin_int_cplx}) and multiplying by the factor of $-1$ in front yields the following equivalent forms to express the result for \cite[Eq.~3.937~3]{gr14}:
\begin{subequations}\label{eq:3.937_3_cmplx}
\begin{align}
&\hspace{-12pt}\int_{0}^{2\pi}  \!\exp(p \cos x{\,+\,}q \sin x) \sin(q \cos x{\,-\,}p \sin x{\,+\,}mx) \, dx \quad \scalebox{0.9}{$[p \textrm{ and } q \textrm{ complex-valued}]$}\hspace{-10pt} \nonumber \\
&\quad = \frac{i \pi}{m!} \Big(\big[(p_R{+}q_I){\,+\,}i(p_I{-}q_R) \big]^m - \big[(p_R{-}q_I){\,+\,}i(p_I{+}q_R) \big]^m \Big) \label{eq:3.937_3_cmplx_v1}\\
&\quad = \frac{i \pi}{m!} \Big(\big[(p_R{+}q_I)^2{\,+\,}(p_I{-}q_R)^2 \big]^{m/2} e^{i\,m \arg[(p_R+q_I)+i(p_I-q_R)]} \label{eq:3.937_3_cmplx_v2} \\[-5pt]
&\hspace{54pt} - \big[(p_R{-}q_I)^2{\,+\,}(p_I{+}q_R)^2 \big]^{m/2} e^{i\,m \arg[(p_R-q_I)+i(p_I+q_R)]} \Big)  \nonumber \\
&\quad = \frac{i \pi}{m!} \big[(p-iq)^m - (p+iq)^m \big]. \label{eq:3.937_3_cmplx_v3}
\end{align}
\end{subequations}
The form in (\ref{eq:3.937_3_cmplx_v3}) is essentially the same as at the start of (\ref{eq:3.937_3_deriv}). However, in this case, since the imaginary parts of $p$ and $q$ may not be zero, the terms cannot be combined to form a $\sin$ term in the same way.

In a similar way, we can insert the constants into (\ref{eq:cos_int_cplx}) to yield the following equivalent forms to express the result for \cite[Eq.~3.937~4]{gr14}:\pagebreak[0]
\begin{subequations}\label{eq:3.937_4_cmplx} 
\begin{align}
&\hspace{-12pt}\int_{0}^{2\pi}  \!\exp(p \cos x{\,+\,}q \sin x) \cos(q \cos x{\,-\,}p \sin x{\,+\,}mx) \, dx \quad \scalebox{0.9}{$[p \textrm{ and } q \textrm{ complex-valued}]$} \hspace{-10pt} \nonumber \\
&\quad = \frac{\pi}{m!} \Big(\big[(p_R{+}q_I){\,+\,}i(p_I{-}q_R) \big]^m + \big[(p_R{-}q_I){\,+\,}i(p_I{+}q_R) \big]^m \Big) \label{eq:3.937_4_cmplx_v1}\\
&\quad = \frac{\pi}{m!} \Big(\big[(p_R{+}q_I)^2{\,+\,}(p_I{-}q_R)^2 \big]^{m/2} e^{i\,m \arg[(p_R+q_I)+i(p_I-q_R)]} \label{eq:3.937_4_cmplx_v2} \\[-5pt]
&\hspace{54pt} + \big[(p_R{-}q_I)^2{\,+\,}(p_I{+}q_R)^2 \big]^{m/2} e^{i\,m \arg[(p_R-q_I)+i(p_I+q_R)]} \Big)  \nonumber \\
&\quad = \frac{\pi}{m!} \big[(p-iq)^m + (p+iq)^m \big] \label{eq:3.937_4_cmplx_v3}
\end{align}
\end{subequations}

\section{Conclusion}\label{sec:concl}
In this paper, we have examined several integrals that appear in Gradshteyn and Ryzhik \cite{gr14}. Our main focus has been on \cite[Eq.~3.937~1]{gr14} \cite[Eq.~3.937~2]{gr14}; the others are special cases of these two integrals. We have determined that the formulas for these two integrals produce a sign error about half the time if the integer $m$ is odd, and have derived the conditions for the formulas' parameters that lead to a sign error. We furthermore have derived updated expressions ((\ref{eq:sin_int_new}) and (\ref{eq:cos_int_new})) that correct the errors, are simpler, and can be used for a wider range of parameter values. For the special cases, we have determined that some are correct but can be generalized further, while others contain errors as well. To summarize:
\begin{itemize}
    \item \cite[Eq.~3.931~4]{gr14}: The minus sign may be replaced by $\pm$.
    \item \cite[Eqs.~3.932~1 and 3.931~2]{gr14}: The formulas are only correct if $m>0$; they are incorrect if $m=0$. \cite[Eq.~3.932~1]{gr14} instead yields 0 if $m=0$, while \cite[Eq.~3.931~2]{gr14} instead yields $\pi$ if $m=0$.
    \item \cite[Eq.~3.936~1]{gr14}: The formula is correct as given.
    \item \cite[Eqs.~3.936~2 and 3.936~3]{gr14}: The restriction $p>0$ is unnecessary and may be removed.
    \item \cite[Eq.~3.936~4]{gr14}: The integrated equation may be generalized to $e^{p \cos x} \times \sin(p\sin x \linebreak[0]\pm mx)$; the same result of 0 will be obtained.
    \item \cite[Eqs.~3.937~3 and 3.937~4]{gr14}: If $m$ is odd, an error in sign will result from the formulas when $p<0$. To correct this and also allow $p=0$ to be used, $\atan(q/p)$ should be replaced by $\operatorname{atan2}(q,p)$ or $\arg(p+iq)$.
\end{itemize}

Lastly, we have also considered the extended case where the parameters are complex-valued rather than real-valued, and have derived the results for the integrals in this event. In the case of \cite[Eqs.~3.937~1 to 3.937~4]{gr14}, the resulting expressions (respectively (\ref{eq:sin_int_cplx}), (\ref{eq:cos_int_cplx}), (\ref{eq:3.937_3_cmplx}), and (\ref{eq:3.937_4_cmplx})) are somewhat more complicated. However, for the other integrals, it turns out that the same formulas (extended and/or corrected) still apply if $p$ is complex-valued.

\section*{Acknowledgments}
The authors would like to thank Dr.~Ivo Maljevi\'{c} of TELUS Communications for his collaboration and his comments on this work.

This work was supported by funding from TELUS Communications and from the Natural Sciences and Engineering Research Council (NSERC) of Canada.

\bibliographystyle{plain}
\bibliography{articles_RCE090123}
\end{document}